\documentclass[11pt]{article}
\usepackage[usenames]{color}
\usepackage[colorlinks=true,
linkcolor=webgreen,
filecolor=webbrown,
citecolor=webgreen]{hyperref}

\definecolor{webgreen}{rgb}{0,.5,0}
\definecolor{webbrown}{rgb}{.6,0,0}

\usepackage{amssymb,psfig,epsfig,latexsym,graphicx,here}

\newtheorem{theorem}{Theorem}

\newcommand{\eqn}[1]{(\ref{#1})}
\newcommand{\hsp}{\hspace*{\parindent}}

\newcommand{\eeq}{\end{equation}}
\newcommand{\beql}[1]{\begin{equation}\label{#1}}
\newcommand{\bsq}{{\vrule height .9ex width .8ex depth -.1ex }}
\newcommand{\NN}{{\mathbb N}}
\newcommand{\ZZ}{{\mathbb Z}}

\makeatletter
\def\@sect#1#2#3#4#5#6[#7]#8{\ifnum #2>\c@secnumdepth
     \def\@svsec{}\else
     \refstepcounter{#1}\edef\@svsec{\csname the#1\endcsname.\hskip .75em }\fi
     \@tempskipa #5\relax
      \ifdim \@tempskipa>\z@
        \begingroup #6\relax
          \@hangfrom{\hskip #3\relax\@svsec}{\interlinepenalty \@M #8\par}%
        \endgroup
       \csname #1mark\endcsname{#7}\addcontentsline
         {toc}{#1}{\ifnum #2>\c@secnumdepth \else
                      \protect\numberline{\csname the#1\endcsname}\fi
                    #7}\else
        \def\@svsechd{#6\hskip #3\@svsec #8\csname #1mark\endcsname
                      {#7}\addcontentsline
                           {toc}{#1}{\ifnum #2>\c@secnumdepth \else
                             \protect\numberline{\csname the#1\endcsname}\fi
                       #7}}\fi
     \@xsect{#5}}
\def\@begintheorem#1#2{\it \trivlist \item[\hskip \labelsep{\bf #1\ #2.}]}
\makeatother

\begin{document}
\begin{center}

{\large {\bf Sloping Binary Numbers: A New Sequence Related to the Binary Numbers}} \\
\vspace*{+.2in}
\end{center}

\begin{center}
David Applegate, \\
Internet and Network Systems Research Center,
AT\&T Shannon Labs, \\
180 Park Avenue,
Florham Park, NJ 07932--0971, USA \\
(Email: david@research.att.com)
\smallskip

Benoit Cloitre, \\
13 rue Pinaigrier, \\
Tours 3700, FRANCE \\
(Email: abcloitre@wanadoo.fr)
\smallskip

Philippe Del\'{e}ham \\
Lyc\'{e}e Polyvalent des Iles, \\
BP 887,
98820 W\'{e} Lifou,
NEW CALEDONIA \\
(Email: kolotoko@lagoon.nc)
\smallskip

N.~J.~A.~Sloane${}^{(1)}$, \\
Internet and Network Systems Research Center,
AT\&T Shannon Labs, \\
180 Park Avenue,
Florham Park, NJ 07932--0971, USA \\
(Email: njas@research.att.com)

\end{center}

\bigskip

\begin{center}
April 29, 2005.
\end{center}

\smallskip

\begin{center}
{\bf Abstract}
\end{center}

If the list of binary numbers is read by upward-sloping diagonals,
the resulting ``sloping binary numbers''
0, 11, 110, 101, 100, 1111, 1010, $\ldots$ 
(or 0, 3, 6, 5, 4, 15, 10, $\ldots$) have some surprising properties.
We give formulae for the $n$-th term and the $n$-th missing term,
and discuss a number of related sequences.

\vspace{0.8\baselineskip}
${}^{(1)}$ To whom correspondence should be addressed.

\vspace{0.8\baselineskip}
Keywords: binary numbers, integer sequences, permutations of integers

\vspace{0.8\baselineskip}
AMS 2000 Classification: Primary 11B83, secondary 11A99, 11B37.

\section{Introduction}
\hsp
We start by writing the binary expansions of the numbers 0, 1, 2, $\ldots$ in an array:
\setlength{\arraycolsep}{2pt}
\renewcommand{\arraystretch}{0.8}
$$\begin{array}{cccc}
&&&0 \\
&&&1 \\
&&1&0 \\
&&1&1 \\
&1&0&0 \\
&1&0&1 \\
&1&1&0 \\
&1&1&1 \\
1&0&0&0 \\
1&0&0&1 \\
&\cdot &\cdot &\cdot
\end{array}
$$
\renewcommand{\arraystretch}{1.0}%
\setlength{\arraycolsep}{5pt}%
By reading this array along diagonals that
slope upwards to the right we obtain the sequence
$$
0, \,  11, \,  110, \,  101, \,  100, \,  1111, \,  1010, \,  1001, \,  1000, \,  1011, \,  \ldots
$$
of {\em sloping binary numbers}, which we denote by $s(0)$, $s(1) , \ldots$.
Written in base 10, $s(0)$, $s(1)$, $s(2), \ldots$ are
$$
0, 3, 6, 5, 4, 15, 10, 9, 8, 11, \ldots \qquad
{\rm (A102370)}\footnote{Six-digit numbers
prefixed by `A' indicate the corresponding entry in \cite{OEIS}.}.
$$
Our goal is to study those numbers as well as several related sequences.
Table \ref{T1} shows $s(0), \ldots, s(32)$ both
in binary and decimal, together with the corresponding values of $(s(n)-n)/2$.
Not every nonnegative number occurs as an $s(n)$ value:
in particular, the numbers 1, 2, 7, 12, 29, 62, 123, 248, 505, $\ldots$
(A102371) never appear.
We denote the omitted numbers by $t(1), t(2), t(3), \ldots$.

In Section \ref{sec2} we state our main theorems,
which give formulae and recurrences for $s(n)$ and $t(n)$,
as well as for a downward-sloping version $d(n)$.
In Section \ref{sec3} we discuss some further properties of these numbers,
namely the trajectories under repeated application of the map $n \mapsto s(n)$ 
(it is interesting that the
trajectory of 2, for example, follows a simple
rule for at least the first 400 million terms, but eventually this rule breaks down);
the fixed points (numbers $n$ such that $s(n) =n$);
the number of terms in the summations in \eqn{EqSF1}
and \eqn{EqTF1} (two number-theoretic functions that
may be of independent interest);
and the average order of $s(n)$.
In the final section, Section \ref{sec4}, we give
two related sequences $\sigma (n)$ and $\delta (n)$
which are permutations of the nonnegative integers,
and a second downward-sloping sequence which is obtained 
by left-adjusting the array of binary numbers.

It is worth mentioning that this work has given rise to
an unusually large number of new sequences---see
the list at the end of this paper.
Only the most important of these will be mentioned in the paper.
Conversely, we were surprised to find very few points of contact with sequences already present in \cite{OEIS}, sequence A034797 being one of the few exceptions.

\section{ The main theorems }\label{sec2}

The first theorem gives the basic properties of the sloping binary numbers $s(n)$.

\begin{theorem}\label{ThS}
\begin{itemize}

\item[$(i)$]
Let $n > 0$. Then for any $m > \log_2 n$,
\beql{EqSF0}
s(n) ~ = ~ 2^m - \frac{1}{2} - \frac{1}{2}
\sum_{k=0}^m
(-1)^{\left\lfloor \frac{n+k}{2^k} \right\rfloor} 2^k ~.
\eeq

\item[$(ii)$]
$s(n)$ satisfies the recurrence
$s(0) =0$ and, for $i \ge 0$, $0 \le j \le 2^i -1$,
\beql{EqSFR}
s(2^i +j) = \left\{
\begin{array}{ll}
2^i + s(j)  & \mbox{~if~} j \ne 2^i - i-1 \\ [+.1in]
3 \cdot 2^i + s(j) & \mbox{~if~} j = 2^i - i-1 ~.
\end{array}
\right.
\eeq

\item[$(iii)$]
\beql{EqSF1}
s(n) ~=~ n ~+~ \sum_{k \ge 1 , \atop n+k \, \equiv \, 0\pmod{2^k}} 2^k \,.
\eeq

\item[$(iv)$]
The values of $s(n)$ are distinct, and $s(n) \ge n$ for all $n \ge 0$.

\end{itemize}
\end{theorem}

\begin{table}[htb]
\caption{The sloping binary numbers $s(n)$ are 
obtained by reading the array of binary numbers along upward-sloping diagonals.
The table gives $s(0), \ldots, s(32)$ in both base 2 and base 10,
as well as the values of $(s(n)-n)/2$.}

\setlength{\arraycolsep}{2pt}
\renewcommand{\arraystretch}{0.8}
$$
\begin{array}{|  @{~} l@{~~~~}llllll|llllll@{~~~~}r @{~} |c|} \hline
\multicolumn{7}{|c|}{n\rule{0in}{.17in} } & \multicolumn{7}{c|}{s(n)} &~(s(n)-n)/2~ \\ [0.04in]
\hline  
0\rule{0in}{.17in} & &&&&&0 \, & \, &&&&&0 & 0 & 0 \\
1* & &&&&&1 \, & \,  &&&&1 & 1 & 3 & 1 \\
2* & &&&&1&0 \, & \, &&&1 &1 & 0 & 6 & 2 \\
3 & &&&&1&1 \, & \, &&&1&0&1&5 & 1 \\
4 & &&&1&0&0 \, & \, &&&1&0&0&4 & 0 \\
5* & &&&1&0&1 \, & \, &&1&1&1&1&15 & 5 \\
6 & &&&1&1&0 \, & \, &&1&0&1&0&10 & 2 \\
7 & &&&1&1&1 \, & \, &&1&0&0&1&9 & 1 \\
8 & &&1&0&0&0 \, & \, &&1&0&0&0&8 & 0 \\
9 & &&1&0&0&1 \, & \, &&1&0&1&1&11 & 1 \\
10 & &&1&0&1&0 \, & \, &&1&1&1&0&14 & 2 \\
11 & &&1&0&1&1 \, & \, &&1&1&0&1&13 & 1 \\
12* & &&1&1&0&0 \, & \, &1&1&1&0&0&28 & 8 \\ 
13 & &&1&1&0&1 \, & \, &1&0&1&1&1&23 & 5 \\
14 & &&1&1&1&0 \, & \, &1&0&0&1&0&18 & 2 \\
15 & &&1&1&1&1 \, & \, &1&0&0&0&1&17 & 1 \\
16 & &1&0&0&0&0 \, & \, &1&0&0&0&0&16 & 0 \\
17 & &1&0&0&0&1 \, & \, &1&0&0&1&1& 19 & 1 \\
18 & &1&0&0&1&0 \, & \, &1&0&1&1&0&22 & 2 \\
19 & &1&0&0&1&1 \, & \, &1&0&1&0&1&21 & 1 \\
20 & &1&0&1&0&0 \, & \, &1&0&1&0&0&20 & 0 \\
21 & &1&0&1&0&1 \, & \, &1&1&1&1&1&31 & 5 \\
22 & &1&0&1&1&0 \, & \, &1&1&0&1&0&26 & 2 \\
23 & &1&0&1&1&1 \, & \, &1&1&0&0&1&25 & 1 \\
24 & &1&1&0&0&0 \, & \, &1&1&0&0&0&24 & 0 \\
25 & &1&1&0&0&1 \, & \, &1&1&0&1&1&27 & 1 \\
26 & &1&1&0&1&0 \, & \, &1&1&1&1&0&30 & 2 \\
27* & &1&1&0&1&1 \, & \, 1&1&1&1&0&1&61 & 17 \\
28 & &1&1&1&0&0 \, & \, 1&0&1&1&0&0&44 & 8 \\
29 & &1&1&1&0&1 \, & \, 1&0&0&1&1&1&39 & 5 \\
30 & &1&1&1&1&0 \, & \, 1&0&0&0&1&0&34 & 2 \\
31 & &1&1&1&1&1 \, & \, 1&0&0&0&0&1&33 & 1 \\
32 & 1&0&0&0&0&0 \, & \, 1&0&0&0&0&0&32 & 0 \\ [0.03in]
\hline
\end{array}
$$
\label{T1}
\end{table}
\renewcommand{\arraystretch}{1.0}
\setlength{\arraycolsep}{5pt}

\paragraph{Proof.}
We first establish some notation. If the
binary expansion of a nonnegative number $n$ is
$$
n = a_0 + a_1 2 + a_2 2^2 + ... + a_m 2^m ~,
$$
where $a_k \in \{0,1\}$, then we call
$a_k$ the $2^k$'s bit of $n$. For future reference
we note that
\beql{EqBIN1}
a_k = \frac{1-(-1)^{\left\lfloor \frac{n}{2^k} \right\rfloor}}{2}, \quad k \ge 0 ~,
\eeq
and so
\beql{EqBIN2}
n ~=~ \sum_{k = 0}^{\infty} 
\, \frac{1-(-1)^{\left\lfloor \frac{n}{2^k} \right\rfloor}}{2} \, 2^k ~,
\eeq
where the upper limit in the summation
can be replaced by $\lfloor \log _2 n \rfloor$.
In Theorem \ref{ThT} we will use
the $2$'s-complement binary expansion for numbers $n < 0$.
This is obtained by writing the binary expansion 
of the nonnegative number $-(n+1)$ as a string beginning
with infinitely many $0$'s, and replacing all $0$'s by $1$'s
and all $1$'s by $0$'s.
Thus the binary expansion of a negative number begins with infinitely many $1$'s
(see Table \ref{T2} below).

(i) Let $L$ denote the infinite, right-adjusted,
array formed from the binary expansions of the
nonnegative numbers (as on the left of Table \ref{T1}),
and let $R$ be the corresponding array formed by
the binary expansions of $s(0), s(1), \ldots$
(as in the central column of the table).
It follows at once from the definition of $s(n)$ that
the right-hand columns (the $1$'s bits)
of $L$ and $R$ agree, the second column from the right in
$R$ (the $2$'s bits) is obtained by shifting the $2$'s column
of $L$ upwards by one place, the $4$'s column of $R$ is
obtained by shifting the $4$'s column of $L$ upwards by two places,
the $8$'s column by three places, and so on.

We also see from Table \ref{T1} that while there are $2^k$ vectors
$u\in \{0,1\}^k$ in $L$, there are only $2^k -1$ such vectors in $R$.
Exactly one vector $u \in \{0,1\}^k$ is missing from each set of $2^k$:
this is $t(k)$.

Because of the way the columns of L are shifted
to form R,  we have (compare \eqn{EqBIN2}):
\beql{EqBIN3}
s(n) ~=~ \sum_{k = 0}^{\infty} 
\, \frac{1-(-1)^{\left\lfloor \frac{n+k}{2^k} \right\rfloor}}{2} \, 2^k ~,
\eeq
where now the upper limit in the summation
can be replaced by any number $m > \log _2 n$.
Therefore, for such an $m$, we have
$$
s(n) ~=~
\frac{2^{m+1}-1}{2} ~-~
\frac{1}{2} \,
\sum_{k = 0}^{m} 
\, (-1)^{\left\lfloor \frac{n+k}{2^k} \right\rfloor} \, 2^k ~,
$$
which proves \eqn{EqSF0}.

(ii)  We will prove \eqn{EqSFR} for $i \ge 2$,
the cases $i=0$ and $1$ being trivial.
Let $n = 2^i + j$.
We consider three subcases.  \\
(a) If $2^i \le n < 2^{i+1} -(i+1)$, then
the diagonal for $n$ is identical to the diagonal for $j$,
except that the $2^i$'s bit is $1$, so
$s(n) = 2^i + s(j)$. \\
(b) If $ n = 2^{i+1} -(i+1)$, then
the diagonal for $n$ is identical to the diagonal for $j$,
except that the $2^i$'s and $2^{i+1}$'s bits are $1$, so
$s(n) = 2^{i+1} + 2^i + s(j)$. \\
(c) If $ 2^{i+1} -(i+1) < n < 2^{i+1} $, then
the diagonal for $n$ is identical to the diagonal for $j$, except that
it has a 0 in the $2^i$'s bit and a 1 in the $2^{i+1}$'s bit,
whereas the diagonal for $j$ has
a 1 in the $2^i$'s bit and a 0 in the $2^{i+1}$'s bit.
Therefore
$s(n) = 2^{i+1} - 2^i + s(j)$. In each case \eqn{EqSFR} holds.

(iii) The starred values of $n$ in the first column of
Table \ref{T1} indicate where the $2^k$'s bit of $s(n)$ is set
for the first time.
Let $p_k =2^k -k$.
Then the $2^k$'s bit $(k \ge 1 )$ of $s(n)$ is set,
and is the highest bit set, 
precisely for $n \in \{p_k , p_k +1, \ldots, p_{k+1} -1 \}$.

The effect of the upwards shift of the columns of L
can be expressed in another way.
Consider the values $(s(n) - n)/2$ (see the final column of Table \ref{T1}).
Each such term is a sum.
Starting with the empty sum, if $n$ is odd we add 1 to the sum,
if $n$ is in the arithmetic progression 2, 6, 10, 14, $\ldots$ 
we add 2, and in general, for $k \ge 1$, if $n$ is in
the arithmetic progression $p_k + i2^k$ $(i\ge 0)$ we add $2^{k-1}$.
But $n$ is in this arithmetic progression precisely
when $n+k \equiv 0$ $(\bmod~2^k )$.
Thus
$$
\frac{s(n)-n}{2} ~=~
\sum_{k \ge 1 , \atop n+k \equiv 0 ~(\bmod~2^k )} 2^{k-1} \,,
$$
which proves \eqn{EqSF1}.
Equation \eqn{EqSF1} can also be deduced from \eqn{EqSFR},
using induction on $i$.

(iv) Equation \eqn{EqSF1} implies that $s(n) \ge n$.
It remains to show that the values $s(n)$ are distinct.
Suppose $n \ne m$.
Let $2^i$ be the highest power of $2$ which divides
$n-m$.  Then $n-m = 2^i + j 2^{i+1}$, for some integer $j$, and
$$
\Big\lfloor \frac{n+i}{2^i} \Big\rfloor =
\Big\lfloor \frac{m+2^i+j 2^{i+1}+i}{2^i} \Big\rfloor
= \Big\lfloor \frac{m+i}{2^i} \Big\rfloor + 1 + 2j \, .
$$
From \eqn{EqBIN3},
this means that the coefficients of $2^i$ in the binary expansions of
$s(n)$ and $s(m)$ are different, so $s(n) \ne s(m)$.
This completes the proof of (iv) and of the theorem.~~~$\bsq$

\noindent{\bf Remarks.}
1.~~The argument in the final paragraph of the proof shows that 
if $n$ is not congruent to $m$ mod $2^k$, then $s(n)$ is not congruent to
$s(m)$ mod $2^k$ (since $n$ not congruent to $m$ mod $2^k$ means $i \le k$ in
that argument).  
Therefore all $2^k$ congruence classes of $n$ mod $2^k$ correspond to distinct
congruence classes of $s(n)$ mod $2^k$.
That is, $s(n)$ is odd if and only if $n$ is odd,
\setlength{\arraycolsep}{2pt}
\renewcommand{\arraystretch}{0.8}
$$\begin{array}{ccc}
\mbox{$s(n)$ ends in:} & \mbox{~if~and~only~if~} & \mbox{$n$ ends in:} \\ [+.1in]
0~0 &~~~& 0~0 \\
0~1 &~~~& 1~1 \\
1~0 &~~~& 1~0 \\
1~1 &~~~& 0~1
\end{array}
$$
respectively, 
$$\begin{array}{ccc}
\mbox{$s(n)$ ends in:} & \mbox{~if~and~only~if~} & \mbox{$n$ ends in:} \\ [+.1in]
0~0~0 & ~~~ & 0~0~0 \\
0~0~1 & ~~~ & 1~1~1 \\
0~1~0 && 1~1~0 \\
0~1~1 && 0~0~1 \\
1~0~0 && 1~0~0 \\
1~0~1 && 0~1~1 \\
1~1~0 && 0~1~0 \\
1~1~1 && 1~0~1
\end{array}
\renewcommand{\arraystretch}{1.0}
\setlength{\arraycolsep}{5pt}
$$
respectively, and so on.
In other words,
for each $k=1,2, \ldots$, there is a permutation $\pi_k$ of
the $2^k$ binary vectors of length $k$ such that the
binary expansion of $s(n)$ ends in $u \in \{0,1\}^k$ if
and only if the binary expansion of $n$ ends in $\pi_k (u)$.

2.~~For the summation in \eqn{EqSF1}, if $n \ge 3$, we need only
consider values of $k \le \lceil \log_2 n \rceil$.

\vspace*{+.05in}

Before studying the missing numbers $t(n)$, it is
convenient to introduce a 
downward-sloping analogue of $s(n)$. If
we read the array L by downward-sloping diagonals,
we obtain the sequence $d(n)$, $n \ge 0$, with initial values
$0, 1, 0, 11, 10, 1, 100, 111, 110, 101, 0, 1011, \ldots$,
or in base $10$, 
\beql{EqSeqD}
0; \, 1; \, 0, 3, 2; \, 1, 4, 7, 6, 5; \, 0, 11, 10, 9, 12, 15, 14, 13, 8; \, 3, 18, 17, \ldots, \quad \quad {\rm (A105033)}\,.
\eeq

Unlike $s(n)$, $d(n)$ is manifestly not one-to-one.
However, there are several similarities between the two sequences.

\begin{theorem}\label{ThD}
\begin{itemize}

\item[$(i)$]
Let $m = \lfloor \log_2 n \rfloor$. Then for $n>0$,
\beql{EqDF0}
d(n) ~ = ~ 2^m - \frac{1}{2} - \frac{1}{2}
\sum_{k=0}^m
(-1)^{\left\lfloor \frac{n-k}{2^k} \right\rfloor} 2^k ~.
\eeq

\item[$(ii)$]
$d(n)$ satisfies the recurrence
$d(0) =0, d(1)=1$ and, for $i \ge 1$, $-1 \le j \le 2^i -1$,
\beql{EqDFR}
d(2^i +i+j) = \left\{
\begin{array}{ll}
d(i-1)  & \mbox{~if~} j = -1 \\ [+.1in]
2^i + d(i+j) & \mbox{~if~} 0 \le j \le 2^i -1 ~.
\end{array}
\right.
\eeq

\item[$(iii)$]
\beql{EqDF1}
d(n) = n - \sum_{1 \le k \le \log_2 n,  \atop n \equiv k-1 ~(\bmod~2^k )} 2^k \,.
\eeq

\end{itemize}
\end{theorem}

\paragraph{Proof.}
The proof is parallel to that of Theorem \ref{ThS} and we
omit the details.~~~$\bsq$

The recurrence \eqn{EqDFR} shows that the $d(n)$ sequence has a natural
division into blocks, where the indices of the blocks run from 
$2^i+i-1$ to $2^{i+1}+(i+1)-2 ~(i \ge 1)$. The 
blocks are separated by semicolons in \eqn{EqSeqD}.

\vspace*{+.05in}

We can now identify the missing numbers $t(n)$.

\begin{theorem}\label{ThT}
\begin{itemize}
\item[$(i)$]
For $n \ge 0$,
\beql{EqTF0}
t(n+1) ~ = ~ 2^n - \frac{1}{2} +
\frac{1}{2} \sum_{k=0}^n (-1)^{\left\lfloor \frac{n-k}{2^k} \right\rfloor} 2^k \,.
\eeq
\item[$(ii)$]
$t(n)$ satisfies the recurrence $t(1) = 1, t(2)=2$ and,
for $i \ge 1$, $i \le j \le 2^i +i$,
\beql{EqTFR}
t(2^i +j) = \left\{
\begin{array}{ll}
2^{2^i+i} -2^i + t(i)  & \mbox{~if~} j = i \\ [+.1in]
2^{2^i+j} -2^i -2^j + t(j)  & \mbox{~if~} i < j \le 2^i+i ~.
\end{array}
\right.
\eeq
\item[$(iii)$]
For $n \ge 1$,
\beql{EqTF1}
t(n) ~=~ -n ~+~ \sum_{k \ge 1 , \atop n-k \, \equiv \, 0\pmod{2^k}} 2^k \,.
\eeq
\item[$(iv)$]
For $n \ge 1$,
\beql{EqTD1}
t(n) =  2^n - 1 - d(n-1).
\eeq
\item[$(v)$]
If we define $s(n)$ for all $n \ge \ZZ$ by {\rm \eqn{EqSF1}}, we have
\beql{EqC}
t(n) = s(-n) \quad\mbox{for}\quad n \ge 1 \,.
\eeq
\end{itemize}
\end{theorem}

\paragraph{Proof.}
We have arranged these formulae in the same order as those in
Theorems \ref{ThS} and \ref{ThD}.
But it is convenient to prove them
in a different order.
(iii)
Continuing from the proof of
Part (iv) of Theorem \ref{ThS},
we observe that the missing numbers
are missing precisely because $s(n)$ for a starred value of $n$
has the $2^k$ bit set;
that is, the $k$-th missing number is found by
erasing the $2^k$ bit from $s(2^k -k )$, or in other words,
$$
t(k) = s(2^k -k) -2^k \quad \quad (k \ge 1 ) \,,$$
from which \eqn{EqTF1} follows immediately.
In the sum in \eqn{EqTF1}, the largest contribution is always from 
$k=n$. For the remaining summands, $1 \le k \le \lfloor \log_2 n \rfloor$.

(iv) now follows from \eqn{EqDF1} and \eqn{EqTF1}, (ii)
from \eqn{EqDFR} and \eqn{EqTD1}, and (v) from \eqn{EqSF1} and \eqn{EqTF1}.

Note that, from \eqn{EqTD1}, $t(n)$ can be obtained
by taking the binary expansion of $d(n-1)$ (written with no leading zeros)
and exchanging $0$'s and $1$'s. 
This leads to a second way to interpret $s(-n)$.
Let us write the binary expansions of the negative
numbers (using the $2$'s-complement notation)
above the binary expansions of the nonnegative numbers,
as in Table \ref{T2}.  If we define $s(n)$ for {\em all} $n$
by reading along upward-sloping diagonals,
we see that $s(-1), s(-2), s(-3), \ldots$ are 
$1, 10, 111, 1100, \ldots$, or in base $10$,
the numbers $1, 2, 7, 12, ...$.
That these numbers really are the missing 
numbers $t(1), t(2), t(3), \ldots$
follows from the fact that reading the upper half of Table \ref{T2}
along upward-sloping diagonals
is the same as reversing the order of the rows in the upper half
of the table, exchanging $0$'s and $1$'s, and reading downwards.
That is, $s(-n) = 2^n - 1 - d(n-1) = t(n)$, which we know to be true
from (iv) and (v).

(i) Finally,  we obtain \eqn{EqTF0}  by considering how the columns
in the upper half of Table \ref{T2} have been shifted,
just as we obtained \eqn{EqSF0} by considering how
the columns in the lower half
of the table were shifted. 
This completes the proof of the theorem.~~~$\bsq$

\begin{table}[htb]
\caption{By using 2's-complement notation for
the binary expansion of negative numbers, $s(n)$ can be defined for all $n \in \ZZ$.
The values $\{s(n): n \le -1 \}$ are the numbers missing from
$\{s(n) : n \ge 0\}$.}

$$
\begin{array}{| @{~} r@{~~~}ccccc|cccccc@{~~~}r @{~} |} \hline
\multicolumn{6}{|c|}{n\rule{0in}{.15in} } & \multicolumn{7}{c|}{s(n)} \\ [0.04in]
\hline
-6\rule{0in}{.17in} & \cdots  1 & 1 & 0 & 1 & 0 \, & \, 1 & 1 & 1 & 1 & 1 & 0 & 62 \\
-5 & \cdots  1 & 1 & 0 & 1 & 1 \, & \, & 1 & 1 & 1 & 0 & 1 & 29 \\
-4 & \cdots  1 & 1 & 1 & 0 & 0 \, & \, & & 1 & 1 & 0 & 0 & 12 \\
-3 & \cdots  1 & 1 & 1 & 0 & 1 \, & \, & & & 1 & 1 & 1 & 7 \\
-2 & \cdots  1 & 1 & 1 & 1 & 0 \, & \, & & & & 1 & 0 & 2 \\
-1 & \cdots  1 & 1 & 1 & 1 & 1 \, & \, & & & & & 1 & 1 \\
0 & \cdots   0 & 0 & 0 & 0 & 0 \, & \, & & & & & 0 & 0 \\
1 & \cdots   0 & 0 & 0 & 0 & 1 \, & \, & & & & 1 & 1 & 3 \\
2 & \cdots   0 & 0 & 0 & 1 & 0 \, & \, & & & 1 & 1 & 0 & 6 \\
3 & \cdots   0 & 0 & 0 & 1 & 1 \, & \, & & & 1 & 0 & 1 & 5 \\
4 & \cdots   0 & 0 & 1 & 0 & 0 \, & \, & & & 1 & 0 & 0 & 4 \\
5 & \cdots   0 & 0 & 1 & 0 & 1 \, & \, & & 1 & 1 & 1 & 1 & 15 \\
6 & \cdots   0 & 0 & 1 & 1 & 0 \, & \, & & 1 & 0 & 1 & 0 & 10 \\ [0.03in]
\hline
\end{array}
$$
\label{T2}
\end{table}

\vspace*{+.05in}

\noindent{\bf Remarks.}
1.~~Since the values $\{s(-n) = t(n): n \ge 1\}$ are the
numbers missing from the sequence $\{s(n): n \ge 0\}$,
$s$ is a bijection from the integers $\ZZ$ to the nonnegative integers $\NN$.
The inverse map $s^{-1}$ is a bijection from $\NN$ to $\ZZ$, 
with initial values 
$s^{-1} (0)$, $s^{-1} (1)$, $s^{-1} (2) , \ldots$ given by
$$
\begin{array}{l}
0, -1, -2, 1,4,3,2,-3,8,7,6,9,-4,11,10,5,16,15,14,17,20,19, \\[+.1in]
~~~~18,13,24,23,22,25,12,-5,26, \ldots \qquad\qquad {\rm (A103122)}
\end{array}
$$

\noindent
2.~~The periodicity of the columns of Table \ref{T2} shows that
the permutations $\pi_k$ relating the final $k$ bits 
of $n$ and $s(n)$ also relate the final $k$ bits of $n$ and $t(n)$.

\noindent
3.~~It is worth mentioning the coincidence which led
us to discover \eqn{EqTF1}.
We considered the sequence
\beql{EqD}
R(k) ~:=~ s(p_k) ~=~ 2^k -k \, + 
\sum_{l \ge 1 , \atop k \, \equiv \, l \pmod{2^l}} 2^l \,, ~~~ k \ge 1 \
\eeq
(the values of $s(n)$ which exceed a new power of 2, see Table \ref{T1}),
which begins 3, 6, 15, 28, 61, 126, $\ldots$ (A103529).
Both $R(k)$ and $t(k)$ are just less than powers of 2,
and to our surprise it appeared from the numerical data that
\beql{EqE}
2^{k+1} - R(k) = 2^k - t(k) , \quad k \ge 1 \,,
\eeq
taking the values 
\beql{EqE2}
1, 2, 1, 4, 3, 2, 5, 8, 7, 6, \ldots, \quad\quad {\rm (A103530)} \,,
\eeq
and this coincidence (which
is a consequence of Theorem \ref{ThT}) suggested \eqn{EqTF1}.

We end this section by listing some further formulae
relating these numbers.
They follow easily from the above theorems.

\begin{itemize}

\item[(i)]
For $n \ge 0$ and any $j$ with $j \le n < 2^j$,
\beql{EqCross1}
d(n) = 2^j - 1 -s(2^j-1-n) \,.
\eeq

\item[(ii)]
For $n \ge 0$ and any $j$ with $0 \le n < 2^j-j$,
\beql{EqCross2}
s(n) = 2^j - 1 -d(2^j-1-n) \,.
\eeq

\item[(iii)]
\beql{EqCross3}
t(n) = s(2^n - n) - 2^n, ~~n\ge 1 \,.
\eeq

\end{itemize}

\section{Further properties}\label{sec3}
\hsp
In this section we discuss some further properties of these sequences.
\subsection{Trajectories}\label{sec3.1}
\hsp
Let $T_m = \{m, s(m), s(s(m)), s(s(s(m))), \ldots \}$ 
denote the trajectory of $m$ under repeated application of the map $n \mapsto s(n)$.
The initial terms of $T_m$ {\em appear} to follow simple rules.
For example,
$$
T_1 = 1,3,5,15,17,19,21,31,33,35,37,47, \ldots, \quad\quad {\rm (A103192)} \,,
$$
{\em appears} to agree with the sequence $\widehat{T}_1$ of
numbers that are congruent to $-1, 1, 3$ or $5$ mod $16$ (A103127).
In fact these two sequences agree precisely for the first 511 terms:
\renewcommand{\arraystretch}{1.2}
$$
\begin{array}{rrrr}
n & ~T_1 (n) & ~\widehat{T}_1 (n) & \mbox{Difference} \\ [0.03in] \hline
0 & 1 & 1 & 0 \\
1 & 3 & 3 & 0 \\
2 & 5 & 5 & 0 \\
\cdots & \cdots & \cdots & \cdots \\
510 & 2037 & 2037 & 0 \\
511 & 4095 & 2047 & 2048 \\
512 & 4097 & 2049 & 2048 \\
\cdots & \cdots & \cdots & \cdots
\end{array}
$$
\renewcommand{\arraystretch}{1.0}

The explanation for this lies in the following theorem.

\begin{theorem}\label{Th2}
For $n$ in any arithmetic progression $\{aj+b : j \ge 0 \}$,
where $a \ge 1$ and $b$ are integers,
the values $s(n) -n$ are unbounded.
\end{theorem}

\noindent{\bf Proof.}
Let $a = c2^d$ with $c$ odd.
For any $m$ such that $m \ge d$ and $2^m > b+d$,
let $k=2^m -b$ and choose $j \ge 1$ so that
$$
cj \equiv -2^{m-d} ~\bmod~2^{2^m -b-d}
$$
(this has a solution since $c$ is odd).
Then for $n= aj+b$ it is easy to check that $n+k \equiv 0$ (mod $2^k$),
and so $s(n) \ge n+2^k$.~~~$\bsq$

\vspace*{+.05in}
This phenomenon is shown more dramatically in $T_2$, which begins
$$
2,6,10,14,18,22,26,30,34,38,42,46,50,54,58,126,130,134, \ldots \qquad
{\rm (A103747)}
$$
The initial terms match the sequence $\widehat{T}_2$ defined by
$$
\widehat{T}_2 (16j+i) := 8(16j+i) + \epsilon_i ~,
$$
for $j \ge 0$, $0 \le i \le 15$, where $\epsilon_0 , \ldots, \epsilon_{15}$ are
$$
2, -2, -6, -10, -14, -18, -22, -26, -30, -34, -38, -42, -46, -50, -54, 6\,.
$$
We have checked by computer that the sequences $T_2$
and $\widehat{T}_2$ agree for at least 400 million terms.
On the other hand, the above theorem shows that the sequences must eventually diverge.
For suppose on the contrary that $T_2 (n) = \widehat{T}_2 (n)$ for all $n$, 
and consider the arithmetic progression $128j + 2$, $j \ge 0$.
These are the values $\widehat{T}_2 (16j)$, and in $\widehat{T}_2$ are followed by
$128j+6$.
But the proof of Theorem \ref{Th2} shows that 
when $j= 2^{119} -1$, $s(128j+2) \ge 128j +2+2^{126} \neq 128j +6$.
So certainly by term $n = 8(2^{119} -1) \approx 10^{36.72 \ldots}$,
$T_2$ and $\widehat{T}_2$ disagree.

\subsection{Fixed points}\label{sec3.2}
\hsp
The fixed points of $s(n)$, that is,
the numbers $n$ for which $s(n) =n$, are observed to be
\beql{EqFA}
0,4,8,16,20,24,32,36,40,48,52, \ldots, \qquad\quad {\rm (A104235)}
\eeq
Dividing by 4 we obtain
\beql{EqFB}
0,1,2,4,5,6,8,9,10,12,13,14,16, \ldots, \qquad\quad {\rm (A104401)}
\eeq
which omits the numbers
\beql{EqFC}
3,7,11,15,19,23,27,31,35, \ldots, \quad\quad {\rm (A103543)}
\eeq
The latter sequence in fact consists of the numbers
of the form $4j +3$ $(j \ge 0)$, together with
\beql{EqFD}
62, 126, 190, 254, 318, 382, 446, 510, 574, 638, \ldots . \quad\quad {\rm (A103584)}
\eeq

The following theorem explains these observations.

\begin{theorem}\label{ThF}
$s(n) = n$ if and only if $n \equiv 0$ $(\bmod~4)$
and $n$ does not belong to any of the arithmetic progressions
\beql{EqFIX1}
Q_r:= \{ 2^{4r}j - 4r : j \ge 1 \} \, ,
\eeq
for $r = 1, 2, \ldots$.
\end{theorem}

\noindent{\bf Proof.}
These are straightforward verifications using \eqn{EqSF1}, which shows
that $s(n) > n$ if and only if
$n + k \equiv 0$ $(\bmod~2^k)$ 
for some $k \geq 1$.  
From $k=1$ and $k=2$,
we have that if $s(n) = n$ then $n \equiv 0$ $(\bmod~4)$.  
We may exclude $k \geq 3$ with
$k \not\equiv 0$ $(\bmod~4)$
because such $k$ are subsumed by $k=1$ and $k=2$.~~~$\bsq$

\vspace*{+.05in}
\noindent{\bf Remark.}
An examination of \eqn{EqTF1} shows that we may restrict \eqn{EqFIX1}
to $r$ such that $t(4r) = 2^{4r} - 4r$, since if $t(4r) \not= 2^{4r} -
4r$, $Q_r$ will be contained in $Q_s$ for some $s < r$.

\subsection{The number of terms in the formulae for $t(n)$ and $s(n)$}
\hsp
While studying the missing numbers $t(n)$,
we investigated the number $f(n)$
(say) of terms in the summation in \eqn{EqTF1},
or, equally, in the summation \eqn{EqD} for the record values $R(k)$.
That is,
\beql{EqTA}
f(n) ~:=~ \# \{1 \le k \le n: ~k \equiv n \mbox{~mod~} 2^k \, \}\, ,
\eeq
a number-theoretic function which may be of independent interest.
The initial values $f(1), \ldots, f(32)$ are
$$
\begin{array}{cccccccccccccccccc}
1 & 1 & 2 & 1 & 2 & 2 & 2 & 1 & 2 & 2 & 3 & 1 & 2 & 2 & 2 & 1 \\
2 & 2 & 3 & 2 & 2 & 2 & 2 & 1 & 2 & 2 & 3 & 1 & 2 & 2 & 2 & 1 & \cdots & ~{\rm (A103318)}
\end{array}
$$
The smallest $n$ such that $f(n) =3$ is 11,
corresponding to the values $k=1 , 3$ and 11;
and a 4 appears for the first time at $f(2059)$.
But although we computed 20 million terms, we were unable to find a 5.
This is explained by the following:
\begin{theorem}\label{ThTA}
$f(n)$ satisfies the
recurrence $f(1) =1$, and, for $i \ge 0$, $1 \le j \le 2^i$,
\beql{EqTB}
f(2^i +j) = \left\{
\begin{array}{ll}
f(j) +1 & \mbox{~if~}  1 \le j \le i  \\ [+.1in]
f(j) & \mbox{~if~}  i+1 \le j \le 2^i \,.
\end{array}
\right.
\eeq
\end{theorem}

\noindent{\bf Proof.}
Note that $k=n$ is always a solution to $k \equiv n$ $(\bmod~2^k )$,
but there are no other solutions with $k > \log_2 n$.
Suppose first that $1 \le j \le i$.
For values of $k$ in the range $1 \le k \le i$,
the equation $k \equiv 2^i +j$ $(\bmod~2^k )$ is equivalent 
to $k \equiv j$ $(\bmod~2^k)$, giving $f(j)$ solutions.
For $k$ in the range $i +1 \le k \le n$,
we get just one further solution, $k=n$, so $f(2^i +j) = f(j) +1$.
On the other hand, suppose that $i+1 \le j \le n$.
We would get $f(j) +1$ solutions, as in the previous case,
except that some values of $k$ that contribute to $f(j)$ are now lost.
The lost values are those $k > \log_2 (2^i +j )$, that is,
$k \ge i+1$.
There is just one such value, namely $k=j$,
and so $f(2^i +j) = f(j)$, as claimed.~~~$\bsq$

\vspace*{+.1in}
\noindent{\bf Corollary.}
{\em
Let $g(m)$ be the minimal value of $n$ such that $f(n) =m$.
Then $g(m)$ satisfies the recurrence $g(1) =1$,
\beql{EqTC}
g(m+1) = 2^{g(m)} + g(m) , \quad m \ge 1 \,,
\eeq
with values
$$
1, \, 3, \, 11, \, 2059, \, 2^{2059} + 2059, \, 
2^{2^{2059} +2059} + 2^{2059} + 2059, \, \ldots,
\quad \quad {\rm (A034797)} \,.
$$
}

\noindent{\bf Proof.}
This is an easy consequence of Theorem \ref{ThTA}, and we omit the details.
It is helpful to use the recurrence \eqn{EqTB} 
to build up a table of values of $f(n)$, as shown in Table \ref{TG}.
For instance, the values $f(17)$ through $f(32)$
are obtained by copying the values $f(1)$ through $f(16)$,
but adding 1 to the four values
$f(17)$, $f(18)$, $f(19)$, $f(20)$---see
the right-hand column of the central portion of Table \ref{TG}.
The values $f(9)$ through $f(16)$ are likewise obtained by
copying the values $f(1)$ through $f(8)$,
but adding 1 to the three values $f(9)$, $f(10)$, $f(11)$---see
the next-to-last column of the central portion.
In general, the $k$-th column of the central portion, 
starting the count at $k=0$, is, for $k \ge 1$,
periodic with period $2^{k+1}$, the repeating string 
consisting of $2^k$ 0's, $k$ 1's, and then $2^k -k$ 0's.
The total number of 1's in the $n^{\rm th}$ row
of the central part of the table is $f(n)$.
We will not see a row with 4 1's until the $k=11$ column is filled in,
when $f(2^{11} +11 ) = 4$ becomes visible.
To avoid any confusion, we emphasize that the column of
Table \ref{TG} builds up $f(n)$ as a sum of 1's.
It is not the binary expansion of $f(n)$.~~~$\bsq$

\clearpage
\renewcommand{\arraystretch}{0.9}
\begin{table}[htb]
\caption{$f(n)$, the number of 
$1 \le k \le n$ such that $k \equiv n \mbox{~mod~} 2^k$.
The central columns show how $f(n)$ is built up recursively using \eqn{EqTB}.
 }
$$
\begin{array}{| @{~} c @{~} | @{~~} cccccc @{~~} |c|} \hline
n \rule{0in}{.17in} & &&&&&& f(n) \\ [0.03in] \hline
1 \rule{0in}{.16in} & \,  1 & 0 & 0 & 0 & 0 & 0 \, & 1 \\ \hline
2  \rule{0in}{.16in}  & \,  1 & 0 & 0 & 0 & 0 & 0 \, & 1 \\ \hline
3  \rule{0in}{.16in}  & \,  1 & 1 & 0 & 0 & 0 & 0 \, & 2 \\
4 & \,  1 & 0 & 0 & 0 & 0 & 0 \, & 1 \\ \hline
5  \rule{0in}{.16in}  & \,  1 & 0 & 1 & 0 & 0 & 0 \, & 2 \\
6 & \,  1 & 0 & 1 & 0 & 0 & 0 \, & 2 \\
7 & \,  1 & 1 & 0 & 0 & 0 & 0 \, & 2 \\
8 & \,  1 & 0 & 0 & 0 & 0 & 0 \, & 1 \\ \hline
9  \rule{0in}{.16in}  & \,  1 & 0 & 0 & 1 & 0 & 0 \, & 2 \\
10 & \,  1 & 0 & 0 & 1 & 0 & 0 \, & 2 \\
11 & \,  1 & 1 & 0 & 1 & 0 & 0 \, & 3 \\
12 & \,  1 & 0 & 0 & 0 & 0 & 0 \, & 1 \\
13 & \,  1 & 0 & 1 & 0 & 0 & 0 \, & 2 \\
14 & \,  1 & 0 & 1 & 0 & 0 & 0 \, & 2 \\
15 & \,  1 & 1 & 0 & 0 & 0 & 0 \, & 2 \\
16 & \,  1 & 0 & 0 & 0 & 0 & 0 \, & 1 \\ \hline
17   \rule{0in}{.16in} & \,  1 & 0 & 0 & 0 & 1 & 0 \, & 2 \\
18 & \,  1 & 0 & 0 & 0 & 1 & 0 \, & 2 \\
19 & \,  1 & 1 & 0 & 0 & 1 & 0 \, & 3 \\
20 & \,  1 & 0 & 0 & 0 & 1 & 0 \, & 2 \\
21 & \,  1 & 0 & 1 & 0 & 0 & 0 \, & 2 \\
22 & \,  1 & 0 & 1 & 0 & 0 & 0 \, & 2 \\
23 & \,  1 & 1 & 0 & 0 & 0 & 0 \, & 2 \\
24 & \,  1 & 0 & 0 & 0 & 0 & 0 \, & 1 \\
25 & \,  1 & 0 & 0 & 1 & 0 & 0 \, & 2 \\
26 & \,  1 & 0 & 0 & 1 & 0 & 0 \, & 2 \\
27 & \,  1 & 1 & 0 & 1 & 0 & 0 \, & 3 \\
28 & \,  1 & 0 & 0 & 0 & 0 & 0 \, & 1 \\
29 & \,  1 & 0 & 1 & 0 & 0 & 0 \, & 2 \\
30 & \,  1 & 0 & 1 & 0 & 0 & 0 \, & 2 \\
31 & \,  1 & 1 & 0 & 0 & 0 & 0 \, & 2 \\
32 & \,  1 & 0 & 0 & 0 & 0 & 0 \, & 1 \\ [.03in] \hline
33 \rule{0in}{.16in} & \, 1 & 0 & 0 & 0 & 0 & 1 \, & 2 \\
\cdot  & \,\cdot & \cdot &\cdot & \cdot & \cdot & \cdot \, & \cdot \\  [.03in] \hline
\end{array}
$$
\label{TG}
\end{table}
\renewcommand{\arraystretch}{1.0}

\clearpage

\vspace*{+.1in}
\noindent{\bf Remarks.}
1.~~The Corollary explains why our computer search failed to find $g(5)$. \\
2.~~The earliest reference to the sequence $g(m)$ that we have found is 
the entry A034797 in \cite{OEIS}, due to Joseph~L. Shipman,
where it arises as the index of the first impartial game of value $m$,
using the natural enumeration of impartial games (cf. \cite{ONAG}).
It is always rash to make such statements,
especially in view of the connections between games
and coding theory described in \cite{CS86},
but there does not seem to be any connection between 
the present work and the theory of impartial games.

We briefly mention the companion sequence $f' (n)$ (say),
giving the number of terms in the summation in \eqn{EqSF1}.
The initial values $f' (0)$, $f' (1) , \ldots$ are
$$
0,1,1,1,0,2,1,1,0,1,1,1,1,2,1,1,0, \ldots, \quad\quad {\rm (A104234)}\,.
$$
An argument similar to that used to establish Theorem \ref{ThTA} shows:
\begin{theorem}\label{ThTA2}
$f'(n)$ satisfies the recurrence $f' (0) =0$,
$f' (1) =1$, and, for $i \ge 1$, $0 \le j \le 2^i -1$,
\beql{EqTD2}
f' (2^i +j) = \left\{
\begin{array}{ll}
f' (j) +1 & \mbox{~for~} j=2^i - i-1 \\ [+.1in]
f' (j) & \mbox{~otherwise}\,.
\end{array}
\right.
\eeq
Also
\beql{EqTD2a}
f'(2^n-n) \, = \, f(n) \,.
\eeq
\end{theorem}

\vspace*{+.05in}
The positions $g' (0)$, $g' (1), \ldots$
where $f' (n) =0,1,2,3, \ldots$ for the first time are
$$
0,1,5,2037, \ldots, \quad\quad {\rm (A105035)} \,.
$$
We have not investigated this function, but these four values suggest 
the conjecture that $g' (m) = 2^{g(m)} - g(m)$,
which is consistent with \eqn{EqTD2a}.
If so, this would imply that $g' (4) = 2^{2059} - 2059$.
Certainly $f' (2^{2059} - 2059) = 4$, but is this the earliest occurrence
of $4$?

\subsection{Average order}
As the above discussion of trajectories illustrates, the function
$s(n)$ for $n \ge 0$ is quite irregular.
But it is straightforward to compute its average order (cf. \cite[\S18.2]{HW}).

\begin{theorem}\label{Th3}
The average order of $s(n)$ is $n+O (\log \,n )$.
\end{theorem}

The proof is an easy computation from \eqn{EqSF1}.

\section{Related sequences}\label{sec4}
In this section we describe some related sequences.

\subsection{Two permutations of the nonnegative integers}\label{secPerm}
Returning to the standard array of binary numbers,
as on the left of Table \ref{T1}, we define two sequences
related to $s(n)$ and $d(n)$ which are actually
permutations of the nonnegative integers.

The first sequence, $\sigma (n)$, $n \ge 0$, begins
$$
0; 1; 3, 2 ; 6,5,4,7; 15, 10, 9,8,11,14,13,12; \ldots, \quad\quad {\rm (A105027)}
$$
with a block structure indicated by semicolons.
The initial term is 0.
After that, the $m$-th block $(m \ge 0)$,
$$
\sigma (2^m ) , ~~\sigma (2^m +1) , \ldots, \sigma (2^{m+1} -1 ) \,,
$$
is constructed by starting at the leading 1-bit of the numbers
$2^m , \ldots, 2^{m+1} -1$ and reading diagonally upwards and to the right.
The $m$-th block is in fact equal to the terms
$$
s(p_m ) , ~ s(p_m +1) , \ldots, s(p_{m+1} -1), ~ t(m+1) \,,
$$
which relates it to our sequences $s(n)$ and $t(n)$.
For example, the third block consists of the numbers
15, 10, 9, 8, 11, 14, 13, 12, ending with $t(4) = 12$.

On the other hand, if instead we read downwards and to the right,
we obtain the sequence $\delta(n)$, $n \ge 0$, beginning
$$
0;1;3,2;4,7,6,5; 11, 10, 9, 12, 15, 14, 13, 8; \ldots . \quad\quad {\rm (A105025)}
$$
We recognize this as being obtained from $d(n)$ by omitting
repeated terms (compare \eqn{EqSeqD}).
Both $\{\sigma (n) : n \ge 0\}$ and $\{\delta (n) : n \ge 0\}$
are permutations of the nonnegative integers.

\subsection{A second downward-sloping version}\label{sec4.1}

We might have begun by 
{\em left-adjusting} the array of binary numbers,
so that it looks like
\setlength{\arraycolsep}{2pt}
\renewcommand{\arraystretch}{0.8}
$$
\begin{array}{cccc}
0 \\
1 \\
1 & 0 \\
1 & 1 \\
1 & 0 & 0 \\
1 & 0 & 1 \\
1 & 1 & 0 \\
1 & 1 & 1 \\
1 & 0 & 0 & 0 \\
1 & 0 & 0 & 1 \\
1 & 0 & 1 & 0 \\
1 & 0 & 1 & 1 \\
\multicolumn{4}{c}{\cdots}
\end{array}
$$
\renewcommand{\arraystretch}{1.0}%
\setlength{\arraycolsep}{5pt}%
Now if we read by downward-sloping diagonals, we obtain the sequence
$0,10, 110, 101, \ldots$, or in decimal,
$$
0,2,6,5,4,14,13,8,11,10, 9, 12, 30, \ldots \quad\quad
{\rm (A105029)}\,.
$$
This seems less interesting than the previous sequence,
and we have not analyzed it in detail.
There are no repetitions, and the numbers $2^m -1$, $m \ge 1$, do not appear.

\vspace*{+.05in}

Further related sequences can be found 
in the list appended to the end of this paper.

\bigskip
\rule{6in}{.2mm}
\bigskip
\\
\noindent
[Related sequences:
A034797,
A102370,
A102371,
A103122,
A103127,
A103185,
A103192,
A103202,
A103205,
A103318,
A103528,
A103529,
A103530,
A103542,
A103543,
A103581,
A103582,
A103583,
A103584,
A103585,
A103586,
A103587,
A103588,
A103589,
A103615,
A103621,
A103745,
A103747,
A103813,
A103842,
A103863,
A104234,
A104235,
A104378,
A104401,
A104403,
A104489,
A104490,
A104853,
A104893,
A105023,
A105024,
A105025,
A105026,
A105027,
A105028,
A105029,
A105030,
A105031,
A105032,
A105033,
A105034,
A105035,
A105085,
A105104,
A105108,
A105109,
A105153,
A105154,
A105158,
A105159,
A105228,
A105229,
A105271.]

\end{document}